\documentclass [12pt,a4paper,reqno]{amsart}
\usepackage{amsmath,amsthm}
\usepackage{amsfonts}
\usepackage{amssymb}
\usepackage{tikz-cd}
\allowdisplaybreaks

\newcommand{\be}{\begin{equation}}
\newcommand{\ef}{\end{equation}}
\chardef\bslash=`\\ 

\newtheorem*{thm*}{Theorem}

\theoremstyle{definition}

\newtheorem*{remark*}{Remarks}
\newtheorem*{defn*}{Definition}

\theoremstyle{remark}

\newcommand{\wt}{\widetilde}
\newcommand{\wh}{\widehat}
\newcommand{\fc}{\frac}

\newcommand{\iy}{\infty}
 \renewcommand{\sectionmark}[1]{}
\renewcommand{\Im}{\operatorname{Im}}

\newcommand{\ve}{\varepsilon}

\newcommand{\uTs} {universal Teichm\"{u}ller space}
\newcommand{\const}{\operatorname{const}}

 \usepackage{amsfonts}
\newcommand{\field}[1]{\mathbb{#1}}

\newcommand{\D}{\field{D}}
\newcommand{\om}{\omega}
\newcommand{\z}{\zeta}
\newcommand{\ov}{\overline}
\newcommand{\vp}{\varphi}
\newcommand{\hC}{\wh{\field{C}}}
\newcommand{\C}{\field{C}}
\newcommand{\R}{\field{R}}

\newcommand{\B}{\mathbf{B}}
\newcommand{\T}{\mathbf{T}}
\newcommand{\Belt}{\operatorname{Belt}}
\newcommand{\Hol}{\operatorname{Hol}}

\newcommand{\Teich}{\operatorname{Teich}}

\newcommand{\vk} {\varkappa}
\newcommand{\x} {\mathbf x}
\renewcommand{\a} {\alpha}

\begin{document}

\title{Grunsky operator, Grinshpan's conjecture and universal Teichm\"{u}ller space}
\author{Samuel L. Krushkal}

\begin{abstract} A. Grinshpan posed a deep conjecture on the norm of the Grunsky operator generated by univalent functions in the disk. It gives a quantitative answer in terms of the Grunsky coefficients, to
which extent a univalent function determines the bound of dilatations of its quasiconformal
extensions.
We provide the proof of this conjecture and its various analytic, geometric and potential applications.

Another result concerns the model of universal Teichm\"{u}ller space by Grunsky coefficients.
\end{abstract}

\date{\today\hskip4mm({GriCon.tex})}

\maketitle

\bigskip

{\small {\textbf {2020 Mathematics Subject Classification:} Primary: 30C55, 30C62, 30F60, 32F45, 32G15; Secondary 46G20}

\medskip

}

\bigskip

\markboth{Samuel L. Krushkal}{Grunsky operator, Grinshpan's conjecture and universal Teichm\"{u}ller space} \pagestyle{headings}

\bigskip\bigskip
\centerline{\bf 1. PRELIMINARIES}

\bigskip\noindent
{\bf 1.1.  The Grunsky operator on univalent functions}.
We consider the class $S_Q(\iy)$ of univalent functions $f(z) = z + a_2 z^ + \dots$ in
the unit disk $\D = \{|z| < 1\}$ admitting quasiconformal extensions to the whole Riemann sphere
$\hC = \C \cup \{\iy\}$ and its completion $S(\iy)$ in the topology of locally uniform convergence in $\D$.

The Beltrami coefficients of extensions are supported in the complementary disk
$$
\D^* = \hC \setminus \ov{\D} = \{z \in \hC: \ |z| > 1\}
$$
and run over the unit ball
$$
\Belt(\D^*)_1 = \{\mu \in L_\iy(\C): \ \mu(z)|D = 0, \ \ \|\mu\|_\iy  < 1\}.
$$
Each $\mu \in \Belt(\D^*)_1$ determines a unique homeomorphic solution to the Beltrami equation
$\ov \partial w = \mu \partial w$ on $\C$ (quasiconformal automorphism of $\hC$)  normalized by
$w^\mu(0) = 0, \ (w^\mu)^\prime(0) = 1, \ w^\mu(\iy) = \iy$, whose restriction to $\D$
belongs to $S_Q(\iy)$.

The {\bf Schwarzian derivatives} of these functions
$$
S_f(z) = \Bigl(\fc{f^{\prime\prime}(z)}{f^\prime(z)}\Bigr)^\prime
- \fc{1}{2} \Bigl(\fc{f^{\prime\prime}(z)}{f^\prime(z)}\Bigr)^2, \quad f(z) = w^\mu(z)|\D,
$$
belong to the complex Banach space $\B = \B(\D)$ of hyperbolically bounded holomorphic functions in
the disk $\D$ with norm
$$
\|\vp\|_\B = \sup_\D \ (1 - |z|^2)^2 |\vp(z)|
$$
and run over a  bounded domain in $\B$ modeling the {\bf universal Teichm\"{u}ller space} $\T$.
Its origin (the base point) $\vp = \mathbf 0$ corresponds to the identity map $f(z) \equiv z$.
The space $\B$ is dual to the Bergman space $A_1(\D)$, a subspace of $L_1(\D)$ formed by integrable holomorphic functions (quadratic differentials $\vp(z) dz^2$ on $\D$.

One defines for any $f \in S_Q(\iy)$ its {\bf Grunsky coefficients} $\a_{m n}$ from the expansion
 \be\label{1}
\log \fc{f(z) - f(\z)}{z - \z} =  \sum\limits_{m, n = 1}^\iy \a_{m n} z^m \z^n, \quad (z, \z) \in \D^2,
\end{equation}
where the principal branch of the logarithmic function is chosen. These coefficients satisfy
the inequality
 \be\label{2}
\Big\vert \sum\limits_{m,n=1}^\iy \ \sqrt{m n} \ \a_{m n}(f) x_m x_n \Big\vert \le 1
\end{equation}
for any sequence $\mathbf x = (x_n)$ from  the unit sphere $S(l^2)$ of the Hilbert space $l^2$
with norm $\|\x\| = \bigl(\sum\limits_1^\iy |x_n|^2\bigr)^{1/2}$; conversely, the inequality (1)
also is sufficient for univalence of a locally univalent function in $\D$ (cf. \cite{Gr}).

The minimum $k(f)$ of dilatations $k(f^\mu) = \|\mu\|_\iy$ among all quasiconformal extensions
$w^\mu(z)$ of $f$ onto the whole plane $\hC$ (forming the equivalence class of $f$) is called the \textbf{Teichm\"{u}ller norm} of this function. Hence,
$$
k(f) = \tanh d_\T(\mathbf 0, S_f),
$$
where $d_\T$ denotes the Teichm\"{u}ller-Kobayashi metric on $\T$.
This quantity dominates the \textbf{Grunsky norm}
$$
\vk(f) = \sup \Big\{\Big\vert \sum\limits_{m,n=1}^\iy \ \sqrt{mn} \
\a_{m n} x_m x_n \Big\vert: \ \mathbf x = (x_n) \in S(l^2) \Big\}
$$
by $\vk(f) \le k(f)$ (see, e.g., \cite{Kr10}, \cite{Ku2}).
These norms coincide only when any extremal Beltrami coefficient
$\mu_0$ for $f$ (i.e., with $\|\mu_0\|_\iy = k(f)$) satisfies
 \be\label{3}
\|\mu_0\|_\iy = \sup \ \Big\{ \Big\vert \iint_{\D^*} \mu(z) \psi(z) dx dy \Big\vert: \ \psi \in A_1^2(\D^*), \ \|\psi\|_{A_1} = 1\Big\} = \vk(f) \quad (z = x + i y).
\end{equation}
Here $A_1(\D^*)$ denotes the subspace in $L_1(\D^*)$ formed by integrable holomorphic functions
(quadratic differentials) on $\D^*$ (hence, $\psi(z) = c_4 z^{-4} + c_5 z^{-5} + \dots$), so $\psi(z) = O(z^{-4})$ as $z \to \iy$, and $A_1^2(\D^*)$ is its subset consisting
of $\psi$ with zeros even order in $\D^*$, i.e., of the squares of holomorphic functions.
Due to \cite{Kr3}, every $\psi \in A_1^2(\D^*)$ has the form
$$
\psi(z) = \fc{1}{\pi} \sum\limits_{m+n = 4}^\iy \sqrt{m n} \ x_m x_n z^{-(m+n)}
$$
and $\|\psi\|_{A_1(\D^*)} = \|\x\|_{l^2} = 1, \ \x = (x_n)$.

In the general case, $\mu_0 \in \Belt(\D^*)_1$ is extremal in its class if and only if
$$
\|\mu_0\|_\iy = \sup \ \Big\{ \Big\vert \iint_{D^*} \mu(z) \psi(z) dx dy \Big\vert: \ \psi \in A_1(\D^*), \ \|\psi\|_{A_1} = 1\Big\}.
$$
Moreover, if $\vk(f) = k(f)$ and the equivalence class of $f$ is a {\bf Strebel point} of $\T$,
which means that this class contains the Teichm\"{u}ller extremal extension  $f^{k|\psi_0|/\psi_0}$
with $\psi_0 \in A_1(\D)$, then necessarily
 $\psi_0 = \om ^2 \in A_1^2$ (cf. \cite{Kr6}, \cite{Kr10}, \cite{Ku3},\cite{St}).

An important fact is that {\it the Strebel points are dense in any Teichm\"{u}ller space} (see \cite{GL}).

\bigskip
Every Grunsky coefficient $\a_{m n}(f)$ in (1) is represented as a polynomial of a finite number of the initial Taylor coefficients
$a_2, \dots, a_s$ and hence depends holomorphically on Beltrami coefficients
$\mu_f \in \Belt(\D^*)_1$ and on the Schwarzians $S_f \in \T$.
This generates holomorphic maps
 \be\label{4}
 h_\x(S_f) = \sum\limits_{m,n=1}^\iy \ \sqrt{m n} \ \a_{m n} (S_f) \ x_m x_n : \ \T \to \D
\end{equation}
with fixed $\x = (x_n) \in l^2$ with $\|\x\| = 1$ so that
$\sup_{\x\in S(l^2)} h_\x(S_f) = \vk(f)$.

The holomorphy of these functions follows from the holomorphy of coefficients $\a_{m n}$ with respect to Beltrami coefficients $\mu \in \Belt(D^*)_1$ and from the estimate
$$
\Big\vert \sum\limits_{m=j}^M \sum\limits_{n=l}^N \ \beta_{mn} x_m
x_n \Big\vert^2 \le \sum\limits_{m=j}^M |x_m|^2
\sum\limits_{n=l}^N |x_n|^2,
$$
which holds for any finite $M, N$ and $1 \le j \le M, \ 1 \le l \le N$ (see \cite{Po1}, p. 61).

Due to \cite{KK2}, \cite{Kr10}, the set of $f$ with $\vk(f) < k(f)$  is dense in $\Sigma$,
and moreover, the Schwarzian derivatives of these functions form an open and dense set in the universal Teichm\"{u}ller space $\T$.
On the other hand, the functions with $\vk(f) = k(f)$ play a crucial role in applications of Grunsky inequalities to the Teichm\"{u}ller space theory.

Both norms $\vk(f)$ and $k(f)$ are continuous plurisubharmonic functions of $S_f$ on the space $\T$
(in $\B$ norm); see, e.g. \cite{Kr10}.

Note that the Grunsky (matrix) operator
$$
\mathcal G(f) = (\sqrt{m n} \ \a_{mn}(f))_{m,n=1}^\iy
$$
acts as a linear operator $l^2 \to l^2$ contracting the norms of elements $\x \in l^2$;
the norm of this operator equals $\vk(f)$ (cf. \cite{Gr}).

\bigskip
A deep theorem of Pommerenke and Zhuravlev states that any univalent function $f \in S$, with $\vk(f) \le \kappa < 1$ has $\kappa_1$-quasiconformal extensions to $\hC$ with $\kappa_1 = \kappa_1(\kappa) \ge \kappa$ (see \cite{Po1}; [KK1, pp. 82-84), \cite{Zh}). Some explicit estimates for $k_1(k)$ are established in \cite{Kr7}, \cite{Ku5}.

The method of Grunsky inequalities was generalized in several directions, even to bordered Riemann surfaces $X$ with a finite number of boundary components. The case of arbitrary quasidisks is described in \cite{Kr10}.

In particular, for univalent functions $f(z) = z + a_2 z^2 + \dots$ in an arbitrary disk
$\D_r = \{|z| < r\}, \ 0 < r < \iy$, the corresponding function (1) on $(z, \z) \in \D_r^2$
provides the Grunsky norm
$$
\vk(f) = \sup \Big\{\Big\vert \sum\limits_{m,n=1}^\iy \ \sqrt{mn} \ \a_{m n} r^{-m-n} x_m x_n \Big\vert: \ \mathbf x = (x_n) \in S(l^2) \Big\},
$$
and accordingly, the holomorphic maps
$$
 h_{\x,r}(S_f) = \sum\limits_{m,n=1}^\iy \ \sqrt{m n} \ \a_{m n} (S_f) r^{-m-n} \ x_m x_n : \ \T \to \D
$$
with $\sup_{\x\in S(l^2)} h_{\x,r}(S_f) = \vk(f)$.

\bigskip\noindent
{\bf 1.2. The root transform}.
One can apply to $f \in S_Q(\iy)$ the rotational conjugation
$$
\mathcal R_p: \ f(z) \mapsto f_p(z) := f(z^p)^{1/p} = z + \fc{a_2}{p} z^{p-1} + \dots
$$
with integer $p \ge 2$, which transforms $f \in S(\iy)$ into $p$-symmetric univalent functions accordingly
to the commutative diagram
$$
\begin{tikzcd}
\wt \C_p \arrow{d}{\pi_p}  \arrow{r}{\mathcal R_p f}
         &\wt \C_p \arrow{d}{\pi_p}      \\
\hC \arrow{r}{f}
         &\hC
\end{tikzcd}
$$
where $\wt \C_p$ denotes the $p$-sheeted sphere $\hC$ branched over $0$ and $\iy$, and the projection $\pi_p(z) = z^p$.

This transform acts on $\mu \in \Belt(\D^*)_1$ and $\psi \in L_1(\D^*)$ by
$$
\mathcal R_p^* \mu = \mu(z^p) \ov z^{p-1}/z^{p-1}, \quad \mathcal
R_p^* \psi = \psi(z^{-p}) p^2 z^{2p-2};
$$
then
$$
k(\mathcal R_p f) = k(f), \quad  \vk(\mathcal R_p f) \ge \vk(f).
$$
It was established by Grinshpan in \cite{G1} that the root transform does not decrease the
Grunsky norm, i.e.
$$
\vk_p(f) := \vk(f_p) \ge \vk(f_p)
$$
(this also follows from the K\"{u}hnau-Schiffer theorem on reciprocity of the Grunsky norm to
the least positive Fredholm eigenvalue of the curve $L = f(|z| = 1)$; see \cite{Ku2}, \cite{Sc}.

Note that the sequence $\vk_p(f), \ p = 2, 3, \dots$, is not necessarily nondecreasing. For example, for the function $f(z) = z/(1 + t z)^2$ with $|t| \le 1$, we have
$$
f_p(t) = z/(1 + t z^p)^{2/p},
$$
and $\vk(f_p) = |t| = k(f_p)$ for even $p$, while $\vk(f_p)< |t| = k(f_p)$ for any odd $p \ge 3$ (see, e.g. \cite{Kr3}).

\bigskip
The Grunsky coefficients of every function $\mathcal R_p f$ also are polynomials of
$a_2, \dots \ , a_l$, which implies, similar to (4), the holomorphy of maps
 \be\label{5}
h_{\x,p}(\mu) = \sum\limits_{m,n=1}^\iy \ \sqrt{m n} \ \a_{m n} (\mathcal R_p f^\mu) \  x_m x_n : \
\Belt(D^*)_1 \to \D
\end{equation}
for any fixed $p$ and any $\x = (x_n) \in S(l^2)$, and
$\sup_{\x\in S(l^2)} h_{\x,p}(\mu) = \vk(\mathcal R_p f^\mu)$.

Every function $h_{\x,p}(\mu)$ descends to a holomorphic functions on the space $\T$,
which implies that the Grunsky norms $\vk(\mathcal R_p f^\mu)$ are continuous and plurisubharmonic on $\T$ \cite{Kr8}.

\bigskip\noindent
{\bf 1.3. Additional remarks on universal Teichm\"{u}ller space}.
The universal Teichm\"{u}ller space $\T = \Teich (\D)$ is the space of
quasisymmetric homeomorphisms of the unit circle $\mathbb S^1$ factorized by M\"{o}bius maps. Its topology and real geometry is determined by  Teichm\"{u}ller metric which naturally arises from extensions of those $h$ to the unit disk.

This space also admits the complex structure of a complex Banach manifold.
This structure is defined on $\T$ by factorization of the ball $\Belt(\D^*)_1$,
letting $\mu_1, \mu_2 \in \Belt(\D)_1$ be equivalent if the corresponding quasiconformal maps
$w^{\mu_1}, w^{\mu_2}$ coincide on the unit circle $\mathbb S^1 = \partial \D$
(hence, on $\ov{\D}$). Such $\mu$ and the corresponding maps $w^\mu$ are called $\T$-{\bf equivalent}.
The equivalence classes $[w^\mu]_\T$ are in one-to-one correspondence with the Schwarzian derivatives $S_{w^\mu}(z), \ \ z \in \D$.
The factorizing projection
$$
\phi_\T(\mu) = S_{w^\mu}: \ \Belt(\D^*)_1 \to \T
$$
is a holomorphic map from $L_\iy(\D^*)$ to $\B$. This map is a split submersion, which means that $\phi_\T$ has local holomorphic sections (see, e.g., [GL]).

The basic intrinsic metric on the space $\T$ is its {\bf Teichm\"{u}ller} metric
$$
\tau_\T (\phi_\T (\mu), \phi_\T (\nu)) = \frac{1}{2} \inf
\bigl\{ \log K \bigl( w^{\mu_*} \circ \bigl(w^{\nu_*} \bigr)^{-1} \bigr) : \
\mu_* \in \phi_\T(\mu), \nu_* \in \phi_\T(\nu) \bigr\};
$$
it is generated by the canonical Finsler structure on $\wt \T$ (in fact on the tangent bundle
$\mathcal T(\T) = \T \times \B$ of $\T$).

The {\bf Carath\'{e}odory} and {\bf Kobayashi} metrics on $\T$ are, as usually, the smallest
and the largest semi-metrics $d$ on $\T$, which are contracted by holomorphic
maps $h: \ \D \to \T$.
Denote these metrics by $c_\T$ and $d_\T$, respectively; then
$$
c_\T(\psi_1, \psi_2) = \sup \{d_\D(h(\psi_1), h(\psi_2)) : \
h \in \Hol(\T, \D)\},
$$
while $d_\T(\psi_1, \psi_2)$ is the largest pseudometric $d$ on $\T$ satisfying
$$
d(\psi_1, \psi_2) \le \inf \{d_\D(0, t) : \ h(0) = \psi_1, \ \text{and} \
h(t) = \psi_2 \ \ h \in \Hol(\D, \T)\},
$$
where $d_\D$ is the hyperbolic Poincar\'{e} metric on $\D$ of Gaussian curvature $- 4$.

\bigskip\noindent
{\bf 1.4. Truncation}. Fix $0 < \rho < 1$ and consider for $\mu \in \Belt(\D^*)_1$ the maps
$$
f_\rho^{\wt \mu}(z) = \rho^{-1} f^\mu(\rho z), \quad z \in \C
$$
with Beltrami coefficients $\wt \mu(z) = \mu(\rho z)$. Truncating these coefficients by
 \be\label{6}
\mu_\rho(z) = \begin{cases}  \mu(\rho z),  \ \ &|z| > 1,   \\
                             0,     & |z| < 1,
\end{cases}
\end{equation}
one obtains a linear (hence holomorphic) map
$$
\iota_\rho: \ \mu \mapsto \mu_\rho: \ \Belt(\D^*)_1 \to  \Belt(\D^*_{1/\rho})_1.
$$
We compare this map with the holomorphic homotopy
$f_t(z) = \fc{1}{t} f(tz)$ with $|t| \le 1$,
which determines for $|t| < 1$ a holomorphic map of the space $\T$ into itself by
$$
S_f(z) = S_{f_t}(z) = t^2 S_f(t z).
$$

This is obtained by applying, for example, the following lemma of Earle \cite{Ea}.

\bigskip\noindent
{\bf Lemma 1}. {\it Let $E, T$ be open subsets of complex Banach spaces $X, Y$ and $B(E)$ be a Banach space of holomorphic functions on $E$ with sup norm. If $\phi(x, t)$ is a bounded map $E \times T \to B(E)$ such that $t \mapsto \phi(x, t)$ is holomorphic for each $x \in E$, then the map $\phi$ is holomorphic.}

\bigskip
The following important lemma is a special case of a more general result from \cite{Kr1} (also presented in \cite{Kr2}, p. 179). It concerns quasiconformal homeomorphisms with $L_\iy$ bounded but integrally small dilatations

Consider in the space $L_p(\C)$ with $p> 2$ the well-known integral operators
$$
T\rho(z) = - \fc{1}{\pi} \iint\limits_\C \fc{\rho(\z) d \xi d\eta}{\z - z}, \quad
\Pi \rho(\z) = - \fc{1}{\pi} \iint\limits_\C \fc{\rho(\z) d \xi d\eta}{(\z - z)^2}
= \partial_z T\rho(z)
$$
assuming for that $\rho$ has a compact support in $\C$. Then the second integral
exists as a Cauchy principal value, and the derivative $\partial_z T$ generically is understanding as distributional. One of the basic fact in the theorey of quasiconformal maps
is that every quasiconformal automorphism $w^\mu$ with $\|\mu\|_\iy = k < 1$ of the Riemann sphere $\hC = \C \bigcup \{\iy\}$  with $\|\mu\|_\iy = k < 1$ is represented in the form
$$
w^\mu(z) = z + T\rho(z),
$$
where $\rho$ is the solution in $L_p$ (for $2 < p < p_0(k)$)  of the integral equation
$$
\rho = \mu + \mu \Pi \rho,
$$
given by the series
$$
\rho = \mu + \mu \Pi \mu + \mu \Pi \mu(\mu \Pi(\mu)) + \dots.
$$
Let $B_{p,R}$ denote the space of functions $f(z)$ on the disk $\D_R, \ f(0) = 0$, with norm
$$
\|f\|_{B_{p,R}} = \sup_{z_1,z_2 \in \ov{\D_R}}  \ \fc{|f(z_1) - f(z_2)|}{|z_1 - z_2|^{1-2/p}}
+ \|\partial_z f\|_{L_p} + \|\partial_{\ov z} f\|_{L_p}.
$$

\bigskip\noindent
{\bf Lemma 2}. {\it Let $f^\mu$ be a quasiconformal automorphism of $\hC$ conformal in the
disk $\D_R^* = \{|z| > R\}$ normalized by $f^\mu(z) = z + \const + O(1/z)$ as $z \to \iy$ and $f^\mu(0) = 0$. Suppose that $\mu$ satisfies
$$
\|\mu\|_\iy = k < 1, \quad \|\mu\|_{L_r} < \ve,
$$
where $r \ge p_0 p(p_0 - p) = r_0(k, p)$ with $p, p_0 > 2$ indicated above, and $\ve$ is small.
Then $f^\mu$ is represented in the form
$$
f^\mu(z) = z - \fc{1}{\pi} \iint\limits_{\D_R} \mu(\z) \Bigl(\fc{1} {\z - z} - \fc{1}{\z}\Bigr)
d \xi d \eta + \om(z, \mu),
$$
where
$$
\|\om\|_{B_p(\ov{\B_R})} \le M(k, p, R) \ve^2,
$$
and the constant $M(k, p, R)$ depends only on $k, p, R$. }

\bigskip
In other words, the bounded Beltrami coefficients, which are integrally small, define quasiconformal variations of the same form, and their remainder terms are uniformly small remainders of order $\ve^2$.

Since the Beltrami coefficients of maps $f^{\mu_r}$ and $f_r$ differ only on the annulus
$\{1 < |z| < 1/r\}$, one can compare these maps by Lemma 2 when $r = |t|$ is close to $1$.

Assume that $f(z)$ is {\bf asymptotically conformal} on the unit circle $\mathbf S^1$; in other words, for any pair of points $a, b$ of the curve $L = f(\mathbf S^1)$,
$$
\max\limits_{z \in L} \frac{|a - z| + |z - b|}{|a - b|} \to 1 \quad \text{as} \ \ |a - b| \to 0,
$$
where $z$ lies between $a$ and $b$.
Such curves are quasicircles without corners, but can be rather pathological (see, e.g., \cite{Po2}, p.249). In particular, all $C^1$-smooth curves are asymptotically conformal.

Then the Schwarzian deraivative has the growth
$S_f(z) = o((1- |z|)^2) \quad \text{as} \ \ |z| \to 1$
(so the function $(1 - |z|^2)^{-2} S_f(z)$ remains continuous under crossing the unit circle
$\mathbf S^1$, hence on $\ov \D$).
This implies that for any small $\ve > 0$ there exists $r_0 = r_0(\ve) < 1$ such that
$$
(1 - |z|^2)^{-2} |S_{f_r}(z)| < \ve \quad \text{for all} \ \ r \in (r_0, 1).
$$
We fix such $r_0$ and take $r > r_0$ satisfying $(1 + r_0)/2 < r < 1$. For such $r$,
we have by Lemma 2 the uniform bound
 \be\label{7}
\|S_{f^{\mu_r}} - S_{f_r}\|_\B = \a_1(1 - r),
\end{equation}
where $\a_1(1 - r) \to 0$ as $r \to 1$.

Combining this with the continuity of Grunsky norm on $\T$ and holomorphy of functions (5) on
$S_f, \ S_{f_{r_0}}$ and $S_{f^{\mu_r}}$, one obtains the estimate
 \be\label{8}
\vk_p(f_r) - \vk_p(f^{\mu_r}) = \a_{2,p}(1 - r) \to 0 \quad \text{as} \ \ r \to 1; \quad p = 1, 2, \dots
\end{equation}
The remainder depends on $p$ and also is uniform under the indicated bounds for $r$ and fixed $r_0$.

Note also that each $\vk_p(f_t)$ is a radial subharmonic function of $t$ on the unit disk $\{|t| < 1\}$;
hence, all $\vk_p(f_t) = \vk_p(f_{|t|})$ are continuous and monotone increasing on $[0, 1]$ and therefore,
$$
\vk_p(f_r) \le \vk_p(f), \quad \lim_{r\to 1} \vk_p(f_r) = \vk_p(f).
$$
All homotopy functions $f_t(z)$ have for $|t| < 1$ the extremal extensions of Teichm\"{u}ller type;
the uniqueness of these extensions of $f_t$ implies that $\lim_{r\to 1} k(f_r) = k(f)$.

\bigskip\noindent
{\bf 1.5. Quasireflections}. The {\bf quasiconformal reflections} (or quasireflections) are the orientation reversing quasiconformal homeomorphisms of the sphere $\hC$ which preserve point-wise some (oriented) quasicircle $L \subset \hC$ and interchange its interior and exterior domains.  which we denote by $D_L$ and $D_L^*$, respectively.

One defines for $L$ its {\bf reflection coefficient}
$$
 q_L = \inf k(f) = \inf \ \| \partial_z f/\partial_{\ov z} f \|_\iy,
$$
taking the infimum over all quasireflections across $L$, and {\bf quasiconformal dilatation}
$$
Q_L = (1 + q_L)/(1 - q_L) \ge 1.
$$
Due to \cite{Ah2}, \cite{Ku4}, this dilatation is equal to the quantity
$Q_L = (1 + k_L)/(1 - k_L)^2$,
where $k_L$ is the minimal dilatation among all orientation preserving quasiconformal automorphisms $f_{*}$ of $\hC$ carrying the unit circle onto $L$, and
$k(f_{*}) = \|\partial_{\ov z} f_{*}/\partial_z f_{*}\|_\iy$.
On the properties of quasireflections and obtained results see, e.g., \cite{Ah2}, \cite{Kr6}, \cite{Ku4}.

\bigskip\bigskip
\centerline{\bf 2. THE GRINSHPAN CONJECTURE. MAIN THEOREM}

\bigskip\noindent
{\bf 2.1 The Grinshpan conjecture}.
A. Grinshpan posed in \cite{G2} the following deep

\bigskip\noindent
\textbf{Conjecture}. {\it For every $f \in S(\iy)$, we have the equality}
$$
\limsup\limits_{p\to \iy} \vk_p(f) = k(f).
$$

Though this conjecture arose from the theory of univalent functions with quasiconformal extensions, it intrinsically relates to many important problems of geometric complex analysis and of Teichm\"{u}ller space theory.

First of all, geometrically this conjecture implies the equality of the Carath\'{e}odory and Kobayashi metrics on the universal Teichm\"{u}ller space, which leads to many important analytic, geometric and potential features of this space as well as of the space of univalent functions.

Somewhat modified (weakened) version of this conjecture is proven in \cite{Kr8}. We shall use here some
results from this paper.

\bigskip\noindent
{\bf 2.2. Main theorem}.
The aim of this paper is to prove the following theorem giving with its corollaries the answers to above problems.

\bigskip\noindent
{\bf Theorem 1}. {\it Every univalent function $f(z) \in S(\iy)$ with Grunsky norm $\vk(f) < 1$ admits quasiconformal extensions $f^\mu$ to the whole sphere $\hC$ with dilatations
$$
 k(f^\mu) \ge \wh \vk(f) := \lim\limits_{m\to \iy} \vk_{2m}(f),
$$
and
\be\label{9}
k(f) = \wh \vk(f) = \limsup \limits_{p \to \iy} \vk_p(f).
\end{equation} }

\bigskip
The quantity $\wh \vk(f)$ can be regarded as the \textbf{limit Grunsky norm} of $f$.

\bigskip
The fact that any $f \in S$ with $\vk(f) < 1$ belongs to $S_Q$ follows from the
Pommerenke-Zhuravlev theorem mentioned above. Theorem 1 strengthens this theorem giving
explicitly the extremal dilatation of admissible quasiconformal extensions (the Teichm\"{u}ller norm of $f$) and proves positively the Grinshpan conjecture.
This theorem also has many other important consequences related to the complex and potential geometries of the universal Teichm\"{u}ller space and to the well-known Ahlfors question on the intrinsic characterization of conformal maps of the disk onto domains with quasiconformal boundaries. The results are
briefly presented in Section 5.

\bigskip\noindent
{\bf 2.3. Example}. Let us first mention that one cannot replace in the statement of Theorem 1 the assumption $f(z) \in S(\iy)$ by $f(z) \in S_Q$ (i.e., drop the third normalization condition), though it does not appear in the Pommerenke-Zhuravlev theorem. The point is that in this case the root transform
$\mathcal R_p$ {\bf can increase} the Teichm\"{u}ller norm.

For {\bf example}, consider the extremal map $g_r$ in Teichm\"{u}ller's
Verschiebungssatz with minimal dilatation among quasiconformal automorphisms of the unit disk, which are identical on the boundary circle and move the origin into a given point
$- r \in (-1, 0)$. Its Beltrami coefficient $\mu_0 = k |\psi_0|/\psi_0$ is defined by $\psi_0$, which is holomorphic and does not vanish on $\D \setminus \{0\}$ and has simple pole at $0$. This $\psi_0$ is orthogonal to all holomorphic quadratic differentials on $\D$ with respect to pairing
$$
\langle \vp, \psi\rangle = \iint_\D (1 - |z|^2)^2 \vp(z) \ov{\psi(z)} dx dy.
$$
For small $r$, the dilatation $k(g_r) = r/2 + O(r^2)$ (the corresponding formula in \cite{Te}, p. 343,
for extremal dilatation contains an error).

This map $g_r$ extends trivially to a quasiconformal map of $\hC$ by $g_r(z) = z$ for $|z| \ge 1$.
Consider the translated map $f_r(z) = g_r(z) + r$. Its restriction to $\D^*$ has dilatation
$k(f_r) = k(g_r) = 0$.

In contrast, the dilatation of the squared map $f_{r,2} := \mathcal R_2 f_r$ equals $r$, since
the differential $\mathcal R_2^* \psi_0$ is holomorphic on $\D$ and therefore the Beltrami  coefficient $\mathcal R_2^* \mu_0 = \mu_{f_{r,2}}$ is extremal for the boundary values $f_{r,2}|\mathbf S^1$ (note that $f_{r,2}(z) = z + r/(2 z) + \dots$ for $|z| > 1$), and
$$
\vk(f_{r,2}) = k(f_{r,2}) = r + O(r^2), \ r \to 0.
$$
Thus,
$$
\limsup\limits_{p\to \iy} \vk_p (f_r) = \vk_2 (f_r) > k(f_r).
$$
Another example was constructed by R. K\"{u}hnau (private communication).

\bigskip\bigskip
\centerline{\bf 3. PROOF OF THEOREM 1}

\bigskip\noindent
{\bf Step 1: Preliminary lemmas}.

Given a function $f \in S_Q(\iy)$, consider its extremal quasiconformal extension $f^{\mu_0}$ to $\D^*$ with Beltrami coefficient $\mu_0 \in L_\iy(\D^*)$ (hence, $k(f) = \|\mu_0\|_\iy$) and assign to this function the quantity
 \be\label{10}
 \a_{\D^*} = \sup \Big\{ \Big\vert \iint\limits_{\D^*} \mu_0(z) \psi(z) dx dy\Big\vert : \ \psi \in A_1^2(D^*), \ \|\psi\|_{A_1(\D^*)} = 1 \Big\},
\end{equation}

\noindent
{\bf Lemma 3}. \cite{Kr2}, \cite{Kr10} {\it The Grunsky norm  $\vk(f)$ of every function $f \in S_Q(\iy)$ is estimated by its Teichm\"{u}ller norm $k = k(f)$ and the quantity (10) via
  \be\label{11}
\vk(f) \le k \fc{k + \a_{\D^*}(f)}{1 + \a_{\D^*}(f) k},
\end{equation}
and $\vk(f) < k$ unless $\a_{\D^*}(f) = \|\mu_0\|_\iy$.

The last equality occurs if and only if $\vk(f) = k(f)$, and if in addition the equivalence class of $f$ (the collection of maps equal to $f$ on $\partial D$) is a Strebel point, then $\mu_0$ is necessarily of the form
$$
\mu_0(z) = \|\mu_0\|_\iy |\psi_0(z)|/\psi_0(z), \quad \psi_0 \in A_1^2(\D^*).
$$ }

The following lemma plays a crucial role in the proof of Theorem 1.

\noindent\bigskip
{\bf Lemma 4}. {\it Let $f \in S_Q(\iy)$, and let $f^{\mu_0}$ be an extremal extension of
$f$ to $\hC$. Then
 \be\label{12}
k(f^{\mu_0}) = \wt \vk(f^{\mu_0}) := \sup_{\mu\in [f^{\mu_0}]}
\ \limsup\limits_{\rho\to 1} \ \sup_p \sup_{\psi\in A_1^2(\D^*), \|\psi\|_{A_1}=1} \
 \Big\vert \iint_{\D^*} \mathcal R_p^*\mu_{0 \rho}(z) \psi(z) dx dy \Big\vert.
\end{equation}  }

\bigskip
This quantity $ \wt \vk(f)$ cannot be replaced by a smaller lower bound for the dilatations of quasiconformal extensions of $f$; we call it the \textbf{outer limit Grunsky norm} of $f$.

\bigskip\noindent
{\bf Proof}. First assume that $\mu_0$ is of Teichm\"{u}ller form, which means
$$
\mu_0(z) = \|\mu_0\|_\iy |\psi_0(z)|/\psi_0(z),
$$
with $\hC$-holomorphic quadratic differential
 \be\label{13}
\psi_0(z) = c_3 z^{-3} + c_4 z^{-4} + \dots, \quad |z| > 1,
\end{equation}
having at most simple pole at the infinite point.

If $c_3 \ne 0, \ c_4 \ne 0$, then, noting that $\wh \vk(\mathcal R_2 f^{\mu_0}) = \wh \vk(f^{\mu_0})$, one can start with the squared map $\mathcal R_2 f^{\mu_0}$  whose defining quadratic differential is of the form
$$
\mathcal R_2^* \psi_0(z^2) = 4 (c_3 z^{-4} + c_4 z^{-6} + \dots)
$$
and has at $z = \iy$ zero of even order. To avoid a complication of notations, we assume that this holds for $\psi_0$  (hence in (13) $c_3 = 0$).
We only need to consider $\psi_0$ with at least two zeros of odd order.

After applying to $f^{\mu_0}$ the root transform, we get the Teichm\"{u}ller map
$\mathcal R_p^* f= f^{k |\mathcal R_p^* \psi_0|/\mathcal R_p^* \psi_0}$ determined by
quadratic differential
$$
\mathcal R_p^* \psi_0 = \psi_0(z^p) p^2 z^{2p-2} = \wh \psi_0(\z).
$$
Note that if $\psi_0(z)$ has a zero at $z_0 \ne 0, \ z_0 \in \D$, then  $\wh \psi_0(\z) = 0$
at the points $\z = z_0^{1/p}$, and $|z_0^{1/p}| \to 1$ as $p \to \iy$.

Assume that $f^{\mu_0}$ is asymptotically conformal on $\mathbf S^1$ and
fix $\rho_j <1$ arbitrarily close to $1$. Pick $p_j$ so large that all zeros of odd order of
$\mathcal R_{p_j}^* \psi_0$ are placed in the annulus $\{1 < |z| < 1/\rho_j\}$.
Then, taking the truncated Beltrami coefficients $(\mathcal R_{p_j}^* \mu_0)_{\rho_j}$ for
$$
\mathcal R_{p_j}^* \mu_0 = k |\mathcal R_{p_j}^* \psi_0|/\mathcal R_{p_j}^* \psi_0,
$$
vanishing in the disks $\D_{1/\rho_j}$, and applying to these coefficients Lemma 3, one obtains that
on each disk $\D_{1/\rho_j}^*$ the corresponding extremal map $f^{(\mathcal R_{p_j}^* \mu_0)_{\rho_j}}$
is determined by holomorphic quadratic differential with zeros of even order.
Therefore,
   \be\label{14}
\vk(f^{(\mathcal R_{p_j}^* \mu_0)_{\rho_j}}) =
\sup_{(x_n) \in S(l^2)} \Big\vert \sum\limits_{m,n=1}^\iy \ \sqrt{m n} \
\a_{m n}(f^{(\mathcal R_{p_j}^* \mu_0)_{\rho_j}}) \  \rho_j^{m+n} x_m x_n \Big\vert.
\end{equation}

Using this equality, one can find the appropriate sequences $\{\rho_n\} \to 1, \ \{p_n\} \to \iy$
and $\{\psi_n\} \in A_1^2(\D^*)$ with $\|\psi_n\|_{A_1} = 1$ such that in the limit as
$n \to \iy$ the above relations result in the desired equality (12). This proves the assertion
of Lemma 4 for functions $f \in S_Q(\iy)$ admitting Teichm\"{u}ller extensions and such that the curve
$f(\mathbf S^1)$ is asymptotically conformal.

\bigskip
In the case of the generic $f(z)$ having Teichm\"{u}ller extension to $\D^*$, we pass to its homotopy stretching $f_\rho(z) = \fc{1}{r} f(rz)$ with $\rho < 1$ and apply the above arguments giving
the prescribed sequences $\{\rho_n\}, \ \{p_n\}, \ \{\psi_n\}$ giving (13) for each $r$.

In view of properties of norms $\vk_p(f_r)$ and $k(f_r)$ indicated in  Section {\bf 1.3}, one can select the appropriate subsequences $\{\rho_{n_j}\}, \ \{p_{n_j}\}, \ \{\psi_{n_j}\}$ so that
 \be\label{15}
k(f) = \wt \vk(f) = \lim\limits_{\rho_{n_j}\to 1} \ \sup_{p_{n_j}}
\sup_{\psi_{n_j} \in A_1^2(\D^*), \|\psi_{n_j}\|_{A_1}=1} \
 \Big\vert \iint_{\D^*} \mathcal R_p^*\mu_{0 \rho_{n_j}}(z) \psi_{n_j}(z) dx dy \Big\vert,
\end{equation}
which implies the equality (12) for $f$. So, Lemma 4 is valid for all Strebel points, hence, for a dense
subset of $f \in S_Q(\iy)$.

\bigskip
Now consider an arbitrary function $f \in S_Q(\iy)$, and let $\mu_0$ be one of its extremal Beltrami coefficients in $\D^*$ (i.e., with minimal $L_\iy$ norm).

Take a neighborhood $U_0$ of $S_{f^{\mu_0}}$ in $\T$, in which the defining projection
$\phi_\T: \ \Belt(\D^*)_1 \to \T$ has a local holomorphic section $s$, and a sequence of the Strebel points $S_{f^{\mu_n}} \to S_{f^{\mu_0}}$ in $\B$ norm. On this neighborhood, the upper semicontinuous normalization of the function $\wt \vk(f)$ is plurisubharmonic (as a function of the Schwarzians $\vp = S_{f^\mu}$).
In view of subharmonicity, we have
$$
\wt \vk(\vp_0) \ge \limsup \wt \vk(\vp_m) \quad \text{for all sequences} \ \ \vp_m \to \vp_0 = S_{f^{\mu_0}}.
$$
The equalities (15) mean geometrically that the functions
 \be\label{16}
h_{\x,1/\rho_j}(S_{f^{\mathcal R_{p_j}^* \mu_{\rho_j}}}) =   \sum\limits_{m,n=1}^\iy \ \sqrt{m n} \
\a_{m n}(f^{(\mathcal R_{p_j^*} \mu)_{\rho_j}}) \  \rho_j^{m+n} x_m x_n
\end{equation}
form a maximizing sequence for the Carath\'{e}odory distance
$$
c_{\Belt(\D^*)_1}(\mathbf 0, \mu_0) = \tanh^{-1} (\|\mu_0\|_\iy), \quad \mu_0 = s(\vp_0),
$$
on the ball $\Belt(\D^*)_1$ (which coincide with the Kobayashi and Teichm\"{u}ller distances
on this ball).
The continuity of these metrics implies
$$
\wt \vk(f^{\mu_n}) \ge c_{\Belt(\D^*)_1}(\mathbf 0, S_{f^\mu_0}) - \ve_n = \|\mu_0\|_\iy - \ve_n,
\quad \ve_n \to 0,
$$
and therefore, $k(f^{\mu_0}) = \wt \vk(f^{\mu_0})$, and these values coincide with right-hand part of (12).
This completes the proof of Lemma 4.

\bigskip
The following lemma strengthens essentially the estimate (8).

\bigskip\noindent
{\bf Lemma 5.} {\it For any function (5), we have
$$
h_{\x,p}(S_{f^{\mu_r}}) - h_{\x,p}(S_{f_r}) = \a(1 - r),
$$
where the function $\a(t) \to 0$ as $t\to 0$ and depends only on $k(f)$.  }

\bigskip\noindent
{\bf Proof}. Any function $h_{\x,p}$ maps the space $\T$ into the unit disk, and since the Kobayashi
distance does not increase under holomorphic maps, both points $h_{\x,p}(S_{f^{\mu_r}})$ and $h_{\x,p}(S_{f_r})$ must lie in the disk $\{|w| \le k\}, \ k = k(f)$. In view of (7), the hyperbolic distance
between these points satisfies
$$
d_\D(h_{\x,p}(S_{f_r}), h_{\x,p}(S_{f^{\mu_r}})) \le d_\T(S_{f^{\mu_r}}, S_{f_r}) = \a_1(1 - r).
$$
For small $1 - r$, this implies that the Euclidean distance between these points also is estimated by a similar function
$$
|h_{\x,p}(S_{f^{\mu_r}}) - h_{\x,p}(S_{f_r})| = \a(1 - r), \quad \a(r) \to 0 \quad \text{as} \ \ r \to 1,
$$
and this function $\a(t)$ depends only on $k$ and $t$.
This provides the assertion of Lemma 5.

\bigskip\noindent
{\bf Step 2: Equality of the generalized Grunsky norms}.

First assume that $f = f_0$ is holomorphic in the closed unit disk $\ov \D$ (hence, in some disk $\D_d$ of radius $d > 1$, and $f^\prime(z) \ne 0$ in $\D_d$).
Then this function admits the Teichm\"{u}ller extremal extension across any circle $\{|z| = d^\prime\}, \
d^\prime < d$, so the Schwarzian $S_{f_0}$ is a Strebel point of the space $\T$ (in other words, the equivalence class of $f_0$ contains a Teichm\"{u}ller extremal map).

For such $f_0$, the arguments applied in the proof of Lemma 4 provides the sequences $\{r_n\} \to 1, \ \{p_n\} \to \iy, \ \{\x_n\} \subset S(l^2)$ and $\{\psi_n\} \in A_1^2$
defining the extremal extensions $f_0^{(\mathcal R_{p_n}^* \mu_0)_{1/r_n}}$ such that
$$
\vk_{p_n}(f_0^{\mu_{0,r_n}}) \ge \wt \vk(f_0^{\mu_{0,r_n}}) - \ve_n = k(f^{\mu_n}) - \ve_n
$$
with $\ve_n, \ve_n^\prime > 0$ monotone decreasing to zero as $n \to \iy$.
Hence, the corresponding holomorphic functions
$h_n(S_{f_0^{\mu_{0,r_n}}}) : \ \T \to \D$ of type (4) or (16) also satisfy
  \be\label{17}
|h_{\x_n}(S_{f_0^{\mu_{0,r_n}}})| \ge \wt \vk(f_0) - \ve_n = k(f_0)  - \ve_n.
\end{equation}
By  Lemmas 2 and the estimate (7) the functions (17) are approximated  by the corresponding functions $h_n(S_{f_{0,r_n}})$ defining $\wh \vk(f_0)$. Then Lemma 5 implies
$$
|h_n(S_{f_{0,r_n}})| \ge \wh \vk(f_0) - \ve_n^\prime, \quad \ve_n^\prime \to 0,
$$
which results in the limit equalities
  \be\label{18}
\wh \vk(f_0) = \wt \vk(f_0) = k(f_0).
\end{equation}

These relations easily extend to all functions $f \in S_Q(\iy)$ having Teichm\"{u}ller extremal extensions
by preliminary using their homotopy functions $f_r(z)$. By the previous step, the equalities (18) are valid
for all these functions, so $\wh \vk(f_r) = k(f_r) \to k(f)$ as $r \to 1$, which implies the assertion of the theorem for the indicated $f$.

In fact, we have established somewhat more: since the function $S_f \to S_{f_t}$ is holomorphic in $t$
as a map from $\T$ into
$\T \times \D$, the relation (12)  already established for the Strebel points implies the equality
$c_\T(\mathbf 0, \vp) = \tau_\T(\mathbf 0, \vp)$
on all Teichm\"{u}ller disks
$$
\{\phi_\T(t |\psi_0|/\psi_0): \ t \in \D, \ \ \psi_0 \in A_1(\D^*)\} \subset \T.
$$

Finally, consider the remained case, when $S_{f^\mu}$ does not be a Strebel point (hence, generically the curve $f^\mu(\mathbf S^1)$ is fractal). Then one can repeat the arguments from the last step in the proof of Lemma 3
taking an extremal quasiconformal extension $f^{\mu}$ of $f$ and its homotopy functions $f_r^{\mu}(z)$, which obey (18). This $S_{f^\mu}$ is approximated by Strebel points $S_{f^{\mu_n}}$, and in view of the above remark, one can now deal with the Carath\'{e}odory distances $c_\T(\mathbf 0, S_{f^{\mu_n}}$ (i.e.,
consider the corresponding maps (16) as the functions of the Schwarzians $S_f \in \T$).

The desired equalities (18), equivalent to (9), again follow by applying the enveloping function (15) and taking the limit as $r \to 1$.
This completes the proof of Theorem 1.

\bigskip\bigskip
\centerline{\bf 4. SECOND PROOF OF LEMMA 4}

\bigskip
Lemma 4 plays a crucial role in the proof of Theorem 1. Thus we provide also alternate proof
of the last step, which does not involve the properties of invariant distances and of the defining projection $\phi_\T$.

Consider an arbitrary function $f \in S_Q(\iy)$, and let $\mu$ be one of its extremal Beltrami coefficients in $\D^*$ (i.e., with minimal $L_\iy$ norm).

Truncate this $\mu$ by (6) close to $\rho$. Restricting the obtained Beltrami coefficient $\mu_\rho$ to the disk $\D_{1/\rho^\prime}^* = \{|z| > 1/\rho^\prime\}$, one obtains, in view of conformality on the annulus
$$
U_{1/\rho,1/\rho^\prime} = \{1/\rho < |z| < 1/\rho^\prime\},
$$
that the equivalence class of $f^\mu|_{\rho^\prime}$ on the disk $\D_{1/\rho^\prime}$, i.e., among the maps
conformal in the disk $\D_{1/\rho^\prime}$ and with the fixed values on the circle $\{|z| = 1/\rho^\prime\}$  is a Strebel point.
It admits Teichm\"{u}ller extension onto the disk $\D_{1/\rho^\prime}^*$ with Beltrami coefficient
$\mu_{\rho,\rho^\prime}$ satisfying
$$
|\mu_{\rho,\rho^\prime}(z)| = \|\mu_{\rho,\rho^\prime}\|_\iy = \|\mu\|_\iy - o(1) < \|\mu\|_\iy,
$$
where, in view of the properties of extremal Beltrami coefficients, $o(1) = \beta_1(\rho^\prime - \rho) \to 0$ as $\rho^\prime \to \rho$ and both $\rho, \ \rho^\prime$ approach $1$.

The arguments from the second step in the above proof of the lemma applied to $f^{\mu_{r,r^\prime}}$ provide the equality
  \be\label{19}
\vk_p(f^{\mathcal R_p^* \mu_{\rho,\rho^\prime}}) =  k(f^{\mathcal R_p^* \mu_{\rho,\rho^\prime}})
\end{equation}
for all $p \ge p_0(\rho)$ such that $\mathcal R_p^* \mu_{\rho,\rho^\prime}$ have no zeros of odd order
in the disk $\D_{1/\rho^\prime}^*$.

We also have for any fixed natural $p \ge 1$,
$$
\lim\limits_{\rho\to 1} \vk_p(f^{\mu_{\rho,\rho^\prime}}) = \vk_p(f^\mu)
$$
and that for a fixed $\rho$ the functions $\vk_p(f^{\mu_\rho})$ are plurisubharmonic with respect to the Schwarzians $S_{f^\mu_\rho}$ and $S_{f^\mu}$ in $\T$ and bounded by $k(f^\mu)$.
This yields that the upper semicontinuous regularization of the function
 \be\label{20}
\vk_0(f^\mu) = \limsup\limits_{\rho \to 1} \sup_p \vk_p(f^{\mu_{\rho}})
\end{equation}
also is plurisubharmonic on the space $\T$. It admits the circular symmetry on each homotopy disk
 $\{S_{f^\mu_t} \}$ (inherited from the symmetry of $f^\mu|\D$) and thus $\vk_0(f_t^\mu)$ is continuous in $|t|$ on $[0, 1]$.

Applying the relations (19), (20), one derives the desired equality (12), which completes the proof of the lemma.

\bigskip\bigskip
\centerline{\bf 5. APPLICATIONS OF THEOREM 1}

\bigskip
The aim of the section is to illustrate the importance of the limit Grunsky norm, which leads to various sharp bounds for univalent functions with quasiconformal extensions and for the basic curvelinear functionals.

\bigskip\noindent
{\bf 5. 1. Ahlfors' problem}.
In 1963, Ahlfors posed in \cite{Ah1} (and repeated in his book \cite{Ah2}) the following question which gave rise to various investigations of quasiconformal extendibility of univalent functions.

\bigskip
\noindent \textbf{Question.} {\it Let $f$ be a conformal map of the disk (or half-plane) onto
a domain with quasiconformal boundary (quasicircle). How can this map be characterized? }

\bigskip
He conjectured that the characterization should be in analytic properties of the logarithmic derivative $\log f^\prime = f^{\prime\prime}/f^\prime$,
and indeed, many results on quasiconformal extensions of holomorphic maps have been established using $f^{\prime\prime}/f^\prime$ and other invariants (see, e.g., the survey \cite{Kr6} and the references there).

This question relates to another still not completely solved problem in geometric complex analysis:

\bigskip
{\it To what extent does the Riemann mapping function $f$ of a Jordan domain $D \subset \hC$ determine the geometric and conformal invariants (characteristics) of the boundary $\partial D$ and of complementary domain $D^* = \hC \setminus \ov D$?}

\bigskip
Theorem 1 implies a natural qualitative answer to all these questions and shows how the inner features of
conformality completely prescribe the admissible distortion under quasiconformal extensions
of function $f$ and determine the hyperbolic features of the universal Teichm\"{u}ller space.
We present this important consequence from Theorem 1 as

\bigskip\noindent
{\bf Corollary 1}. {\it For any function $f$ mapping conformally the unit disk onto a domain with quasiconformal boundary $L = f(|z| = 1)$, the reflection coefficient $q_L$ of the curve $L$ is determined by the limit Grunsky norm of this function via
 \be\label{21}
 \fc{1 + q_L}{1 - q_L} = \left(\fc{1 + \wh \vk(f)}{1 - \wh \vk(f)}\right)^2.
\end{equation}
Hence the right-hand side of (21) represents the minimal dilatation of quasiconformal
reflections across $L$. }

\bigskip\noindent
{\bf 5.2. Invariant metrics on universal Teichm\"{u}ller space}.
We already mentioned in the proof of Lemma 4 the connection between the function (13) and the invariant distances on $\T$. Here we give their explicit representation generated by the original univalent
functions.

It is elementary that the Carath\'{e}odory, Kobayashi and Teichm\"{u}ller metrics of any Teichm\"{u}ller space $\wt T$ are related by
  \be\label{22}
c_{\wt \T}(\cdot, \cdot) \le d_{\wt \T}(\cdot, \cdot) \le \tau_{\wt\T}(\cdot, \cdot)
\end{equation}
(and similarly for their infinitesimal forms).

The fundamental Royden-Gardiner theorem states that the metrics $d_{\wt \T}$ generated and $\tau_{\wt \T}$ coincide on any space $\wt \T$ (see, e.g., \cite{EKK}, \cite{GL}).
On the other hand, due to \cite{Kr11}, the Carath\'{e}odory metric of the universal Teichm\"{u}ller space
$\T$ also coincides with its Teichm\"{u}ller metric, but this does not hold, for example, for finite dimensional Teichm\"{u}ller spaces of dimension greater than $1$, see \cite{Ga}.
Theorem 1 gives a new proof of this  fact and represents these metrics explicitly in terms of $\wh \vk(f)$.

\bigskip\noindent
{\bf Corollary 2}. {\it The Carath\'{e}odory metric of the universal Teichm\"{u}ller space $\T$ coincides with its Teichm\"{u}ller metric; hence all non-expansive holomorphically invariant metrics on the space $\T$ are equal, in particular, for any two point $\vp_1, \ \vp_2 \in \T$,
In particular, for any pair $(\vp_1, \ \vp_2) \in \T \times \T$,
 \be\label{23}
c_\T(\vp_1, \vp_2) = d_\T(\vp_1, \vp_2) = \tau_\T(\vp_1, \vp_2),
\end{equation}
and similarly for the infinitesimal forms of these metrics.   }

\bigskip\noindent
{\bf Proof}. It follows from Theorem 1 and (19) that for any point $S_{f^\mu} \in \T$, we have the equalities
 \be\label{24}
c_\T(S_{f^\mu}, \mathbf 0) = d_\T(S_{f^\mu}, \mathbf 0) = \tau_\T(S_{f^\mu}, \mathbf 0)
= \tanh^{-1} \wh \vk(f^\mu).
\end{equation}
Now consider two arbitrary points $\vp_1 = S_{f_1}$ and $\vp_2 = S_{f_2}$ in $\T$. Since
the universal Teichm\"{u}ller space is a complex homogeneous domain in $\B$,
this general case is reduced to (24) by moving one of these points to the origin $\vp = \mathbf 0$, applying a right translation of the space $\T$.
Such translations preserve the invariant distances and the Teichm\"{u}ller distance which,
together with (24), yields the equalities (23).

In view of maximality of $d_\T$ and minimality of $c_\T$, all intermediate invariant metrics
also obey (22). This completes the proof of Corollary 2.

\bigskip
This Corollary completely determines the complex geometry of the space $\T$. Note also that
the equalities similar to (24) are also valid for the infinitesimal forms of metrics
$c_\T, \ d_\T, \ \tau_\T$, which provides that all these metrics have holomorphic curvature $- 4$.
We will not go here into details, because the proof involves essentially the results lying out
of the framework of this paper.

\bigskip\noindent
{\bf 5.3. Pluricomplex Green function of universal Teichm\"{u}ller space}.
Corollary 2 also determines the potential features of the space $\T$.
We illustrate this by representation of its {\bf pluricomplex Green function}.

Recall that the pluricomplex Green function $g_D(x, y)$ of a domain $D$ in a complex Banach space $X$
with pole $y$ is defined by
$g_D(x, y) = \sup u_y(x) \quad (x, y \in D)$
followed by the upper semicontinuous regularization
$$
v^*(x) = \lim\limits_{\ve\to 0} \sup_{\|x^\prime - x\|<\ve} v(x^\prime).
$$
The supremum  here is taken over all plurisubharmonic functions
$u_y(x): \ D \to [-\iy, 0)$ such that
$u_y(x) = \log \|x - y\|_X + O(1)$
in a neighborhood of the pole $y$.
Here $\|\cdot\|_X$ denotes the norm on $X$, and the remainder term $O(1)$ is bounded from above.
The Green function $g_D(x, y)$ is a maximal plurisubharmonic  function on
$D \setminus \{y\}$ (unless it is identically $- \iy$).

\bigskip\noindent
{\bf Corollary 3.} {\it For every point $\vp = S_f \in \T$, the pluricomplex Green function $g_\T(\mathbf 0, S_f)$ with pole at the origin of $\T$ is given by
 \be\label{25}
g_\T(\mathbf 0, S_f) = \log \wh \vk(f),
\end{equation}
and for any pair $(\vp, \psi)$ of points in $\T$, we have
$$
g_\T(\vp, \psi) = \log \tanh d_\T(\vp, \psi) = \log \tanh c_\T(\vp, \psi) = \log k(\vp, \psi),
$$
where $k(\vp, \psi)$ denotes the extremal dilatation of quasiconformal maps
determining the Teichm\"{u}ller distance between $\vp$ and $\psi$.   }

\bigskip
The equality (25) is in fact a special case of the general equality
$$
g_D(x, y) = \log \tanh d_D(x, y),
$$
which holds for any hyperbolic Banach domain whose Kobayashi metric $d_D$ is logarithmically plurisubharmonic  (cf. \cite{Di}, \cite{Kl}, \cite{Kr5}).

\bigskip\bigskip
\centerline{\bf 6. MODELING UNIVERSAL TEICHM\"{U}LLER SPACE BY GRUNSKI}
\centerline{\bf COEFFICIENTS}

\bigskip
There are several models of the universal Teichm\"{u}ller space. The most applicable is the Bers
model via the domain $\T$ in the Banach space $\B$ of the Schwarzian derivatives; this model was applied above.

We mention here another model constructed in \cite{Kr12} by applying the Grunsky coefficients of univalent functions in the disk.
In this model, the space $\T$ is represented by a bounded domain in a subspace of $l_\infty$ determined by the Grunsky coefficients.
This domain is biholomorphically equivalent to the Bers domain $\T$.

\bigskip
Consider the univalent nonvanishing functions in the disk $\D^* = \{|z| > 1\}$ with hydrodynamical normalization
 \be\label{26}
f(z) = z + b_1 z^{-1} + \dots: \D^* \to \hC \setminus \{0\}.
\end{equation}
and denote their collection by $\Sigma$. The Grunsky \ coefficients $c_{m n}$ of these functions defined
from the expansion
$$
\log \fc{f(z) - f(\z)}{z - \z} = - \sum\limits_{m, n = 1}^\infty
c_{m n} z^{-m} \z^{-n}, \quad (z, \z) \in (\D^*)^2,
$$
and satisfy (2).
Note that $c_{m n} = c_{n m}$ for all $m, n \ge 1$ and $c_{m 1} = b_m$ for any $m \ge 1$.

Denote by $\Sigma_k$ the subset of $\Sigma$ formed by the functions with $k$-quasiconformal extensions to $\D$, and let $\Sigma^0 = \bigcup_k \Sigma_k$.

Since any $f \in \Sigma^0$ does not vanish in $\D^*$, its
inversion $F_f(z) = 1/f(1/z)$ is holomorphic and univalent in the unit disk $\D$; both functions $f$ and $F_f$ have the same Grunsky coefficients $c_{m n}$.

These coefficients span a $\C$-linear space $\mathcal L^0$ of sequences
$\mathbf c = (c_{m n})$ which satisfy the symmetry relation
$c_{m n} = c_{n m}$ and
$$
|c_{m n}| \le C(\mathbf c)/\sqrt{m n}, \quad
C(\mathbf c) = \const < \infty \quad \text{for all} \ \ m, n \ge 1,
$$
with finite norm
 \be\label{27}
\|\mathbf c\| = \sup_{m,n} \ \sqrt{m n} \ |c_{m n}| +  \sup_{\mathbf x = (x_n) \in S(l^2)} \
\Big\vert \sum\limits_{m,n=1}^\infty \ \sqrt{mn} \ c_{mn} x_m x_n \Big\vert.
\end{equation}
Denote the closure of span $\mathcal L^0$ by $\mathcal L$ and note that the limits of convergent sequences
$\{\mathbf c^{(p)} = (c_{m n}^{(p)})\} \subset \mathcal L$ in the norm (27)
also generate the double series
$$
\sum\limits_{m,n=1}^\infty \ c_{m n} z^{-m} \zeta^{-n}
$$
convergent absolutely in the bidisk $\{(z, \z) \in \hC^2: \ |z| > 1, \ |\z| > 1\}$.

It is established in \cite{Kr12} that {\it the sequences $\mathbf c$  corresponding to functions
$f \in \Sigma^0$ having quasiconformal extensions to the disk $\D$ fill a bounded domain
$\wt \T$ in the indicated Banach space $\mathcal L$ containing the origin, and the correspondence
$$
S_f \leftrightarrow \mathbf c = (c_{m n})
$$
determines a biholomorphism of the domain $\wt \T$ onto the space $\T$.  }

In this model, the Grunsky norm of any function from $S_Q(\iy)$ arises as a canonical part of the Banach norm of its representative $\mathbf c$ in $\wt \T$, and the hyperbolic length of the limit Grunsky norm
of this function is equal by Corollary 2 to each of invariant distances in $\wt \T$ between the point
$\mathbf c$ and the origin. This determines the basic features of both Grunsky norms.

The corresponding holomorphic functions
$$
 h_{\x^0}(\mathbf c) = \sum\limits_{m,n=1}^\iy \ \sqrt{m n} \ c_{m n} \ x_m^0 x_n^0
$$
generating the norm $\vk(f)$ become linear on $\wt \T$, which provide some interesting applications.

We establish here the following property of this domain.

\bigskip\noindent
{\bf Theorem 2}. {\it The domain $\wt \T$ is not starlike (with respect to the origin of $\mathcal L$)}.

\bigskip
The problem on starlikness of Teichm\"{u}ller spaces in Bers' embedding was stated in 1970 in the collection of problems in the book \cite{BK}. This problem still does not be solved completely. Its negative solution for the universal Teichm\"{u}ller space $\T$ was given in \cite{Kr4}. This result also covers other models of $\T$ and later has been extended to finite dimensional Teichm\"{u}ller spaces $\T(0, n)$ of sufficiently large dimensions.
Another proof of non-starlikness of $\T$ is given in \cite{Kr9}. Both proofs explicitly
provide the functions violating this property.

\bigskip\noindent
{\bf Proof}. We apply the same construction as in \cite{Kr9}.
Pick unbounded convex rectilinear polygon $P_n$ with finite vertices $A_1, \dots, A_{n-1}$ and $A_n = \iy$. Denote the exterior angles at $A_j$ by $\pi \a_j$ so that $\pi < \a_j < 2 \pi, \ j = 1, \dots,
n - 1$.
The conformal map $f_n$ of the lower half-plane $H^* = \{z: \ \Im z < 0\}$ onto the complementary polygon $P_n^* = \hC \setminus \ov{P_n}$ is represented by the Schwarz-Christoffel integral
 \be\label{28}
f_n(z) = d_1 \int\limits_0^z (\xi - a_1)^{\alpha_1 - 1} (\xi - a_2)^{\alpha_2 - 1} ... (\xi - a_{n-1})^{\alpha_{n-1} - 1} d \xi + d_0,
\end{equation}
with $a_j = f_n^{-1}(A_j) \in \R$ and complex constants $d_0, d_1$; here $f_n^{-1}(\iy) = \iy$. Its  Schwarzian derivative is given by
$$
S_{f_n}(z) = \mathbf b_{f_n}^\prime(z) - \fc{1}{2}  b_{f_n}^2(z) =
\sum\limits_1^{n-1} \frac{C_j}{(z - a_j)^2} -
\sum\limits_{j,l=1}^{n-1} \frac{C_{jl}}{(z - a_j)(z - a_l)},
$$
where $\mathbf b_f = f^{\prime\prime}/f^\prime$ and
$$
C_j = - (\alpha_j - 1) - (\alpha_j - 1)^2/2 < 0, \ \ C_{jl} =
(\alpha_j - 1)(\alpha_l - 1) > 0.
$$
It defines a point of the universal Teichm\"{u}ller space $\T$ modelled as a bounded domain in the space $\B(H^*)$ of hyperbolically bounded holomorphic functions on $H^*$ with norm
$$
\|\vp\|_{\B(H^*)} = \sup_{H^*} |z - \ov z|^2 |\vp(z)|.
$$

By the Ahlfors-Weill theorem \cite{AW}, every $\vp \in \B(H^*)$ with $\|\vp\|_{\B(H^*)} < 1/2$ is the Schwarzian derivative of a univalent function $f$ in $H^*$, and this function has
quasiconformal extension onto the upper half-plane $H = \{z: \ \Im z > 0\}$ with Beltrami coefficient of the form
 \be\label{29}
\mu_\vp(z) = - 2 y^2 \vp(\ov z), \quad \vp = S_f \ (z = x + i y \in H^*)
\end{equation}
called harmonic.

Denote by $r_0$ the positive root of the equation
$$
\frac{1}{2} \Bigl[\sum\limits_1^{n-1} (\alpha_j - 1)^2 +
\sum\limits_{j,l=1}^{n-1} (\alpha_j - 1)(\alpha_l - 1) \Bigr] r^2 -
\sum\limits_1^{n-1} (\alpha_j - 1) \ r - 2 = 0,
$$
and put
$$
S_{f_n,t} = t b_{f_n}^\prime -  b_{f_n}^2/2, \ t > 0.
$$
Then for appropriate vertices $\a_j$ and appropriate Moebius map $\sigma: \ \D* \to H^*$,  we have the following result from \cite{Kr9}.

\bigskip\noindent
{\bf Lemma 6}. {\it For any convex polygon $P_n$, the Schwarzians $r S_{f_n,r_0}$ with $0 <r < r_0$
create the univalent function $w_r = f_n: H^* \to \C$ whose harmonic Beltrami coefficients
$\nu_r(z) = - (r/2) y^2 S_{f_n,r_0}(\overline z)$ in $H$ are extremal in their equivalence classes, and}
  \be\label{30}
k(f_n \circ \sigma) = \vk(f_n \circ \sigma) = \fc{r}{2} \| S_{f_n,r_0}\|_{\B(H^*)}.
\end{equation}

\bigskip
This lemma yields that any function $w_r$ with $r < r_0$ does not admit extremal quasiconformal extensions of Teichm\"{u}ller type; the extremal extensions have harmonic Beltrami coefficients $\mu_{S_{w_r}}$
given by (29).
Therefore, the Schwarzians $S_{w_r}$ for some $r$ between $r_0$ and $1$ must lie outside of the space $\T$; so this space is not a starlike domain in $\B(H^*)$ and in $\B(\D^*)$.

Consider the corresponding Grunsky coefficients and the interval
$$
I = \{t \mathbf c(f_n \circ \sigma): \ 0 \le t \le 1\}.
$$
Since $S_{f_n \circ \sigma}$ is an inner point of $\T$, it has a neighborhood $U_0$ whose intersection
with the image of the interval $I$ in $\T$ contains a non-degenerate subinterval. Hence,
all $t \mathbf c(f_n \circ \sigma)$ with $t$ sufficiently close to $1$ belong to $\wt \T$ (correspond
to univalent functions), and the same holds by for all $t \in [0, r_0]$, but the interval $I$ does not
lie entirely in $\wt \T$. So this domain is not starlike.

\bigskip
{\small\em{ \leftline{Department of Mathematics, Bar-Ilan
University, 5290002 Ramat-Gan, Israel} \leftline{and
Department of Mathematics, University of Virginia,  Charlottesville, VA 22904-4137, USA}}

\end{document}